%% file: main_Arxiv.tex
\documentclass{article}

\usepackage{amsthm}
\input{package.tex}
\input{commands.tex}

\input{style.tex}

\newcommand{\citet}{\cite}
\newcommand{\citep}{\cite}

\usepackage{PRIMEarxiv}

\usepackage[utf8]{inputenc} 

\title{\input{mytitle}
}
\author{
  Luca Berti\\
  CNRS\\
  Universit\'e de Strasbourg,\\
  Strasbourg, France\\
  \texttt{berti@math.unistra.fr} \\
  \And
  Enrico Facca\\
  Centro Ennio de Giorgi,\\
  Scuola Normale Superiore,\\ 
  Pisa, Italy\\
  \texttt{enrico.facca@sns.it} \\
  \And
  Mario Putti \\
  Department of Mathematics ``Tullio Levi-Civita'' \\
  University of Padova \\
  Padova, Italy\\
  \texttt{putti@math.unipd.it} \\
}

\begin{document}
\maketitle

\begin{abstract}
\input{abstract}
\end{abstract}

\newcommand{\sep}{\and}

\input{paper}

\bibliographystyle{unsrt}  
\bibliography{strings,biblio_abbr,pubblication_ef}

\end{document}

%% file: package.tex
 \usepackage{bbm}

\usepackage{amsmath,amsfonts,amssymb}
\usepackage{cleveref}
\usepackage{mathtools}

\usepackage{placeins}
  
\usepackage{enumerate}

\usepackage{xifthen}
\usepackage{easyvector}
\usepackage{epigraph}
\usepackage{algorithm}

\usepackage{adjustbox}

\usepackage{pgfplots}
\pgfplotsset{compat=newest}
\usetikzlibrary{arrows}
\usetikzlibrary{calc}
\usetikzlibrary{shapes}

\usetikzlibrary{fit,arrows}
\usetikzlibrary{positioning}
\usetikzlibrary{automata}

\usepackage{marginnote}

\usepackage{xpatch}
\makeatletter
\patchcmd{\@mn@margintest}{\@tempswafalse}{\@tempswatrue}{}{}
\patchcmd{\@mn@margintest}{\@tempswafalse}{\@tempswatrue}{}{}
\makeatother

\usepackage{changes}

%% file: commands.tex
\newfont{\NUMBERS}{msbm8 scaled\magstep1}

\newcommand{\REAL}{\mathbb{R}}

\usepackage{ifthen}

\newcommand{\Dim}{d}

\newcommand{\REALNN}{\REAL_{\geq 0}}

\newcommand{\Vect}[2][]
{
  \ifthenelse{\equal{#1}{}}
  {\boldsymbol{#2}}
  {{#2}_{#1}}
}
\newcommand{\Matr}[2][]
{
  \ifthenelse{\equal{#1}{}}
  {\boldsymbol{#2}}
  {{#2}_{#1}}
}

\newcommand{\Span}{\ensuremath{\operatorname{span}}}
\newcommand{\Var}{\ensuremath{\operatorname{var}}}
\newcommand{\Err}{\ensuremath{\operatorname{err}}}

\newcommand{\dx}{\, \mathrm{d}x}

\newcommand{\dv}{\, \mathrm{d}V_{\Metric}}

\newcommand{\Grad}[1][]
{
  \ifthenelse{\equal{#1}{}}
  {{\nabla}}
  {{\nabla_{\!#1}}}
}
\newcommand{\Div}[1][]
{
   \ifthenelse{\equal{#1}{}}
   {{\ensuremath{\operatorname{div}}}}
   {{\ensuremath{\operatorname{div}}_{\!#1}}}
}

\def\XXint#1#2#3{{\setbox0=\hbox{$#1{#2#3}{\int}$ }
\vcenter{\hbox{$#2#3$ }}\kern-.6\wd0}}

\newcommand{\ABS}[1]{\left|#1\right|}
\newcommand{\Dt}[1]{\partial_{t}#1}

\newcommand{\Source}{f^{+}}
\newcommand{\Sink}{f^{-}}

\newcommand{\Supp}{\mbox{supp}}
\newcommand{\SourceDomain}{X}
\newcommand{\SinkDomain}{Y}

\newcommand{\Cost}{c}

\DeclareMathOperator{\Dist}{\operatorname{dist}}

\newcommand{\Plan}{\gamma}
\newcommand{\Plans}{\Pi}

\newcommand{\Transport}{\mu}
\newcommand{\Domain}{\Omega}

\newcommand{\Lyap}{\mathcal{L}}

\newcommand{\LCF}{Lyapunov-candidate functional}

\newcommand{\MK}{MK}

\newcommand{\MKEQS}{MK equations}
\newcommand{\OTP}{OT problem}

\newcommand{\Vof}[2]{[#1]^{#2}}

\newcommand{\Radon}{\mathcal{M}}
\newcommand{\RadonPlus}{\mathcal{M}_{+}}
\newcommand{\Lebesgue}[1]{L^{#1}}
\newcommand{\Lspace}[1]{L^{#1}}
\newcommand{\Cachar}{\mathcal{C}}
\newcommand{\Lip}{\operatorname{Lip}}

\newcommand{\Hilb}[1]{H^{#1}}

\newcommand{\Cont}[1][]{
	\ifthenelse{\equal{#1}{}}
	{\Cachar} 
	{\Cachar^{#1}}
}

\newcommand{\ProjMatSymb}{\mathbb{P}}
\newcommand{\ProjMat}[1][]{
  \ifthenelse{\equal{#1}{}}
  {\ProjMatSymb}
  {\ProjMatSymb_{#1}}
}
\newcommand{\IDSymbol}{\mathbb{I}}
\newcommand{\IDtens}[1][]{
  \ifthenelse{\equal{#1}{}}
  {\IDSymbol}
  {\IDSymbol(#1)}
}
\newcommand{\IDMat}{\IDSymbol}

\newcommand{\PolySymb}{\mathcal{P}}

\newcommand{\PONE}{\ensuremath{\PolySymb_{1}}}

\newcommand{\PZERO}{\ensuremath{\PolySymb_{0}}}

\newcommand{\TolTime}{\tau_{\mbox{{\scriptsize T}}}}

\newcommand{\MPsymb}{h}
\newcommand{\MeshPar}[1][]
{
  \ifthenelse{\equal{#1}{}}
    {\MPsymb}
    {\MPsymb_{#1}}
}

\newcommand{\tstep}{k}
\newcommand{\tstepp}{k+1}
\newcommand{\Deltat}[1][]{\ifthenelse{\equal{#1}{}}{\Delta t_{\tstep}}{\Delta t_{#1}}}
\newcommand{\TdensH}{\Tdens_{\MeshPar}}
\newcommand{\PotH}{\Pot_{\MeshPar}}
\newcommand{\Tsymb}{\mathcal{T}}

\newcommand{\Triang}[1][]
{
  \ifthenelse{\equal{#1}{}}
  {\Tsymb}
  {\Tsymb_{#1}}
}
\newcommand{\TriangH}[1][]
{
  \ifthenelse{\equal{#1}{}}
    {\Tsymb_{\MeshPar}}
    {\Tsymb_{#1}}
}
\newcommand{\Trianghh}{\Triang_{\MeshPar/2}}

\newcommand{\Test}{\phi}

\newcommand{\Ene}{\mathcal{E}}
\newcommand{\Wmass}{\mathcal{M}}
\newcommand{\PotOp}{\mathcal{U}}

\newcommand{\FEM}{FEM}
\newcommand{\VspaceSymb}{\mathcal{V}}
\newcommand{\Vspace}[1][]{
  \ifthenelse{\equal{#1}{}}
  {\VspaceSymb_{\MeshPar}}
  {\VspaceSymb_{\MeshPar}(#1)}
}
\newcommand{\Vdim}{N}
\newcommand{\VbaseSymb}{\varphi}
\newcommand{\Vbase}[1][]{
  \ifthenelse{\equal{#1}{}}
  {\VbaseSymb}
  {\VbaseSymb_{#1}}
}
\newcommand{\WspaceSymb}{\mathcal{W}}
\newcommand{\Wspace}[1][]{
  \ifthenelse{\equal{#1}{}}
  {\WspaceSymb_{\MeshPar}}
  {\WspaceSymb_{\MeshPar}(#1)}
}
\newcommand{\Wdim}{M}
\newcommand{\WbaseSymb}{\psi}
\newcommand{\Wbase}[1][]{
  \ifthenelse{\equal{#1}{}}
  {\WbaseSymb}
  {\WbaseSymb_{#1}}
}

\newcommand{\SurfSymbol}{\Gamma}
\newcommand{\Surf}[1][]
{
  \ifthenelse{\equal{#1}{}}
  {{\SurfSymbol}}
  {{\SurfSymbol}_{\!#1}}
}
\newcommand{\Surfh}{\Surf[\MeshPar]}

\newcommand{\Manifold}{M}        
\newcommand{\Mdim}{n}        
\newcommand{\Chart}[1][]
{
  \ifthenelse{\equal{#1}{}}
  {\phi_{\CutPoint}}
  {\phi_{#1}}
}

\newcommand{\NormalSymbol}{\nu}
\newcommand{\Normal}[1][]{
  \ifthenelse{\equal{#1}{}}
  {{\NormalSymbol}}
  {{\NormalSymbol}_{#1}}
}
\newcommand{\NormalH}{\Normal[\MeshPar]}
\newcommand{\testp}{\psi}
\newcommand{\testx}{\xi}
\newcommand{\uFEMdimension}{\Vdim}
\newcommand{\muFEMdimension}{\Wdim}

\newcommand{\CutPoint}{p}
\newcommand{\Point}{x}

\newcommand{\Sphere}[1]{\mathcal{S}^{#1}}

\newcommand{\SphereDist}[1][]
{
  \ifthenelse{\equal{#1}{}}
  {{r}}
  {{r_{#1}}}
}

\newcommand{\Identification}[1][]
{
  \ifthenelse{\equal{#1}{}}
  {\Psi}
  {\Psi^{#1}}
}

\newcommand{\SphereAngle}[1][]
{
  \ifthenelse{\equal{#1}{}}
  {\varphi}
  {\varphi_{{#1}}}
}

\newcommand{\Tdens}{\mu}
\newcommand{\Pot}{u}
\newcommand{\Vel}{v}

\newcommand{\DMK}{DMK}
\newcommand{\OTD}{OT\  density}
\newcommand{\Opt}[1]{#1^*}
\newcommand{\OptTdens}{\Opt{\Tdens}}
\newcommand{\OptPot}{\Opt{\Pot}}
\newcommand{\OptVel}{\Opt{\Vel}}

\newcommand{\VelH}{\Vel_{\MeshPar}}
\newcommand{\OptVelH}{\Opt{\VelH}}
\newcommand{\OptTdensH}{\Opt{\TdensH}}
\newcommand{\OptPotH}{\Opt{\PotH}}

\newcommand{\OptTime}{\Opt{t}}

\newcommand{\Metric}[1][]
{
  \ifthenelse{\equal{#1}{}}
  {g}
  {g(#1)}
}

\newcommand{\Scal}[2]{\langle #1,#2\rangle}

\newcommand{\Norm}[2][]
{
  \ifthenelse{\equal{#1}{}}
  {{|#2|}}
  {{|#2|}_{#1}}
}
\newcommand{\NORM}[2][]
{
  \ifthenelse{\equal{#1}{}}
  {{\|#2\|}}
  {{\|#2\|}_{#1}}
}

\newcommand{\EMDADMM}{EMDADMM}
\newcommand{\BP}{BP}
\newcommand{\Wass}[1]{W_{#1}}
\newcommand{\ErrWass}{\Err_{\Wass{1}}}

\newcommand{\Curve}{\sigma}
\newcommand{\Tangent}[2][]
{
  \ifthenelse{\equal{#1}{}}
  {T#2}
  {T_{#1}#2}
}

\newcommand{\Cellsymb}{\mathbbm{E}}
\newcommand{\Cell}[1][]{\ifthenelse{\equal{#1}{}}{\Cellsymb}{\Cellsymb_{#1}}}
\newcommand{\SubCellsymb}{\mathbbm{e}}
\newcommand{\SubCell}[1][]{\ifthenelse{\equal{#1}{}}{\SubCellsymb}{\SubCellsymb_{#1}}}
\newcommand{\CellH}[1][]{\ifthenelse{\equal{#1}{}}{\Cellsymb_{\MeshPar}}{\Cellsymb_{\MeshPar,#1}}}

\newcommand{\Ipot}{j}

\newcommand{\NCell}[1][]{
  \ifthenelse{\equal{#1}{}}
    {N_{\Cellsymb}}
    {N_{\Cellsymb({#1})}}
}

\newcommand{\Edgesymb}{\sigma}
\newcommand{\Edge}[1][]{
  \ifthenelse{\equal{#1}{}}
    {\Edgesymb}
    {\Edgesymb_{#1}}
}
\newcommand{\EdgeH}[1][]{
  \ifthenelse{\equal{#1}{}}
    {\Edgesymb_{\MeshPar}}
    {\Edgesymb_{\MeshPar,#1}}
}
\newcommand{\InradiusSymbol}{r}
\newcommand{\Inradius}[1][]
{
  \ifthenelse{\equal{#1}{}}
  {\InradiusSymbol}
  {\InradiusSymbol_{\scriptscriptstyle{#1}}}
}

\newcommand{\NeighSymbol}{\mathcal{N}}
\newcommand{\Neigh}[1][]
{
  \ifthenelse{\equal{#1}{}}
  {\NeighSymbol_{\Point}}
  {\NeighSymbol_{#1}}
}

\newcommand{\Itdcell}{r}
\newcommand{\Itd}{r}

\newcommand{\Itdsubcell}{l}

\newcommand{\Of}[1]{(#1)}

\newcommand{\Gen}[1]{\bar{#1}}
\newcommand{\Perr}{1}

\newcommand{\ADMMOptVelH}{\bar{\Vel}_{\MeshPar}}

\newcommand{\Neig}{{N_e}}

\newcommand{\Fgen}{\phi}

\newcommand{\Kmax}{K_{\mbox{max}}}

%% file: style.tex
\newtheorem{Defin}{Definition}

\newtheorem{Prop}{Proposition}

\newtheorem{Problem}{Problem}

\newtheorem{Remark}{Remark}

\Crefname{equation}{Eq.}{Eqs.}
\crefname{equation}{eq.}{eqs.}
\crefname{theorem}{Theorem}{Theorems}
\crefname{chapter}{Chapter}{Chapters}
\crefname{figure}{Figure}{Figures}
\crefname{Conject}{Conjecture}{Conjectures}
\crefname{Problem}{Problem}{Problems}
\crefname{Prop}{Proposition}{Propositions}        
\crefname{Theo}{Theorem}{Theorems}
\crefname{Lemma}{Lemma}{Lemma}
\crefname{Corollary}{Corollary}{Corollaries}
\crefname{section}{Section}{Sections}

%% file: mytitle.tex
Numerical Solution of the $L^1$-Optimal Transport Problem on Surfaces

%% file: abstract.tex
In this article we study the numerical solution of the
$L^1$-Optimal Transport Problem on 2D
surfaces embedded in $\REAL^3$, via the \DMK\ 
formulation introduced in \cite{Facca-et-al:2018}.
We extend from the Euclidean into the Riemannian setting
the \DMK\ model and conjecture the equivalence 
with the solution Monge-Kantorovich equations, a PDE-based
formulation of the $L^1$-Optimal Transport Problem.
We generalize the numerical method proposed in
\cite{Facca-et-al:2018,Facca-et-al-numeric:2020} to 2D
surfaces embedded in $\REAL^3$ using the Surface Finite
Element Model approach to approximate the Laplace-Beltrami
equation arising from the model.
We test the accuracy and efficiency of the proposed numerical
scheme, comparing our approximate solution with respect to an exact
solution on a 2D sphere. The results show that the numerical scheme
is efficient, robust, and more accurate with respect to other
numerical schemes presented in the literature for the solution of
$L^1$-Optimal Transport Problem on 2D surfaces.

%% file: mykeywords.tex
Monge-Kantorovich equations \sep
   Optimal Transport \sep
   Numerical Solution \sep
   Wasserstein Distance \sep
   Riemannian Manifolds \sep
   Triangulated Surfaces

%% file: paper.tex
\section{Introduction}
\label{sec:into}

The \emph{Optimal Transport Problem} (\OTP), first introduced by Monge
in~\cite{Monge:1781} and reformulated by Kantorovich
in~\cite{Kantorovich:1942}, looks for the optimal strategy to
re-allocate a resource from one place ($\SourceDomain$)
to another ($\SinkDomain$) for a given unit-mass
transportation cost $\Cost(x,y)$ ($x\in\SourceDomain$,
  $y\in\SinkDomain$), typically chosen as $\Cost(x,y)=|x-y|^{p}$.
  This problem captured the attention of a number of authors
  (extensive reviews can be found in~\cite{Ambrosio:2003,
    Villani:2003, Villani:2008, Santambrogio:2015}), and, prompted by
  the fluid dynamic interpretation of the quadratic cost problem
  proposed in~\cite{Benamou-Brenier:2000, Benamou-et-al:2002}, a
  number of mathematical and numerical formulations of the \OTP\ have
  been proposed in different fields of applications
  ~\cite{Arjovsky-et-al:2017, Chen:2017, Delzanno-et-al:2008,
    Metivier-et-al:2016, Solomon-et-al:2015}.
  
For $p=1$, i.e., the cost is linear with respect to the distance,
$\Cost(x,y)=|x-y|$, the problem is named $\Lspace{1}$ Optimal
Transport ($\Lspace{1}$-\OTP), and is exactly the problem originally
posed by Monge~\cite{Monge:1781}. In this case, the non-strict
convexity of the cost makes the mathematical analysis more complicated
and increases the difficulties of its numerical solution.  On the other
hand, this problem, and the related problem of calculating the
Wasserstein-1 distance, seems to be more robust with respect to noise
distances defined with other cost functions (e.g., the
quadratic cost $p=2$)~\cite{Peyre:2019}.  The $\Lspace{1}$ problem can
be reformulated in several equivalent ways (see \cite{Ambrosio:2003}
for a complete overview), which can be exploited in the quest for
efficient numerical solvers.

Various numerical methods exploiting these different formulations of
the ($\Lspace{1}$-\OTP) have been developed in recent
years~\cite{Benamou-et-al:2015, Li-et-al:2017}, but the subject is
still a very active field of research.
Among the most widely used numerical solution methods, we can name the
Sinkhorn approach described in~\cite{Cuturi:2013} and its
variants~\cite{Peyre:2019}. These methods are based on entropy
regularization, i.e., they solve a relaxed version of the
linear programming problem associated to the Kantorovich formulation
of the \OTP, and have been successfully applied in many different
fields, such as computer graphics~\cite{Solomon-et-al:2015} or machine
learning~\cite{Genevay-et-al:2016}. One of their limitations is that
they are not suitable for transported densities that are sums of Dirac
measures~\cite{Genevay-et-al:2016}.

Recently,~\cite{Facca-et-al:2018} introduced a dynamic reformulation of
the Monge-Kantorovich equations (\MKEQS), the PDE-based form of the
($\Lspace{1}$-\OTP) proposed by~\cite{Evans-Gangbo:1999}.  The authors
conjecture that the solution of the \MKEQS\ is exactly the long-time
limit of a system of differential equations, called the Dynamic
Monge-Kantorovich (\DMK), that couples a diffusion PDE with an ODE
describing the dynamics of the diffusion coefficient. Unfortunately, a
complete proof of the conjecture is still missing. However, some
theoretical results reported in~\cite{Facca-et-al:2018,
  Facca-et-al-numeric:2020} and extensive numerical evidence presented
in~\cite{Facca-et-al-numeric:2020} strongly support the claim that
\MKEQS\ can be solved via the \DMK\ approach. More
recently,~\cite{Piazzon-et-al:2019} have shown that the equilibrium of
a slight modification of a relaxed version of the \DMK\ indeed converges
in a weak sense to the solution of the \MKEQS, providing additional
support to the above conjecture. The \DMK\ formulation can be easily
solved numerically using standard methods and provides an efficient
numerical approach not only for the solution of two- and
three-dimensional $\Lspace{1}$-\OTP s, but also for the calculation of
the Wasserstein-1 distance between measures~\cite{Facca-et-al-numeric:2020}.

%
%
The $\Lspace{1}$-\OTP\ can be naturally extended from the Euclidean
setting to a Riemannian manifold $\Manifold\subset\REAL^{\Mdim}$ with
metric $\Metric$ by simply replacing the Euclidean distance with the
geodesic distance as the transport cost~\cite{Feldman-McCann:2001,
  Pratelli:2005}.  However, this extension poses non-trivial numerical
challenges, in primis the computation of the geodesic distance
itself. Recent numerical discretizations of \OTP s\ on discrete
surfaces can be found in~\cite{Solomon-et-al:2014, Solomon-et-al:2015,
  Su-et-al:2015, Lavenant:2019}.

In the present paper, we extend the \DMK\ model proposed
in~\cite{Facca-et-al:2018} to the Riemannian setting.  We ambient the
\DMK\ model on a manifold $\Manifold$ with a Riemannian metric
$\Metric$ and address in particular the case in which either $\Sink$
or $\Source$ or both are continuous with respect to the volume form.
The extension of the \DMK\ to Riemannian manifolds is obtained
via an appropriate substitution of the relevant differential operators
with geometry-based equivalents, and takes inspiration
from the definitions introduced in~\cite{Pratelli:2005}. 
Next, we adapt the numerical methods proposed
in~\cite{Facca-et-al-numeric:2020} to solve the DMK on two-dimensional
surfaces embedded in $\REAL^3$. The DMK equations are discretized in
time by means of explicit Euler time stepping and in space by means of
the Surface Finite Element technique developed
in~\cite{Dziuk-Elliott:2013}.  Mimicing the approach used in the
Euclidean (planar) case, the spatial discretization is obtained by a
combination of piecewise constant and conforming piecewise linear
basis functions on nested grids to satisfy a sort of inf-sup stability
condition. The resulting approach has the same computational
properties and costs of the planar $\REAL^2$ equivalent, and thus
provides a very efficient tool for the numerical solution of \OTP\ on
surfaces.

The convergence properties of the developed algorithm are verified on
a test case defined on the unit sphere for whicha closed form solution
of the \MKEQS\ is developed starting from the results reported
in~\cite{Buttazzo-Stepanov:2003}. We experimentally show the
convergence towards steady state of the proposed dynamics and evaluate
the experimental convergence rates of the spatial approximation,
convergence in the approximation of the Wasserstein-$1$ distance
between two measures is tested and is found experimentally to be more
than quadratic with respect to the mesh parameter $\MeshPar$.

Finally, the results obtained with our proposed approach are compared
with those obtained using the method proposed
in~\cite{Solomon-et-al:2014}. This method, tailored to the calculation
of the $\Lspace{1}$-\OTP, also called Earth's Mover Distance (EMD), is
based on the solution of Beckmann problem, which can be shown to be
equivalent to the $\Lspace{1}$-\OTP~\cite{Santambrogio:2015}. Beckmann
problem tries to minimize the $\Lspace{1}$-norm of the vector-valued
measure $\Vel$ subject to the constraint that its divergence must be
equal to the difference between the transported densities. The method
proposed in~\cite{Solomon-et-al:2014}, which we denote with \EMDADMM,
approximates the minimizer of Beckmann problem by means of the
Alternating Direction Method of Multipliers in combination with a
clever FEM discretization of the divergence constraint.  Our
experimental results show that the \DMK\ approach is more accurate
being characterized by a faster convergence rate both in the
approximation of the minimizer of the $\Lspace{1}$-\OTP\ and in the
computation of the Wasserstein-$1$ distance.

The paper is organized as follows. After a short introduction of the
different formulations of the Euclidean $\Lspace{1}$-\OTP\ of interest
for our developments, we summarize the needed variations when the
ambient space is a Riemannian manifold, with particular emphasis to
the DMK approach. Then, we describe the extension of the numerical
formulation of the \DMK\ from the Euclidean to the Riemannian setting
and how surface-\FEM\ can be effectively developed to obtain the
proposed numerical formulation of the $\Lspace{1}$-\OTP.  Finally, the
numerical tests are presented.

\section{The $\Lspace{1}$ optimal transportation problem}
\label{sec:l1opt-mkeqs} 

\subsection{The Kantorovich problem}

The Kantorovich formulation of the $\Lspace{1}$-\OTP\ in the Euclidean
setting can be given as follows.  Consider a convex and compact subset
$\Omega\subset\REAL^\Dim$ with smooth boundary and let $|x-y|$ be
the Euclidean distance between two points $x$ and $y$ in
$\REAL^{\Dim}$. Denote with $\RadonPlus(\Domain)$ the set of
non-negative Radon measures supported in $\Domain$. Given two finite
measures $\Source,\Sink\in\RadonPlus(\Domain)$ such that
$\Source(\Domain)=\Sink(\Domain)$, find the optimal plan $\Opt{\Plan}$
belonging to the set
\begin{equation}
  \label{eq:otp-plan}
  \Plans(\Source,\Sink,\Domain)=
  \left\{ 
    \begin{aligned}
      &\Plan \in
      \RadonPlus(\Domain \times \Domain) :
      \mbox{for all } A,B \mbox{ Borel sets of } \Domain
      \\
      &\Plan(A,\Domain) = \Source(A) \quad
      \Plan(\Domain,B) = \Sink(B)      
    \end{aligned}
  \right\}
\end{equation}
that solves the following minimization problem:
\begin{equation}
  \label{eq:l1-otp}
  \inf_{\Plan}
  \left\{ 
    \int_{\Domain\times\Domain} |x-y| \mathrm{d}\Plan(x,y)
    \ : \
    \Plan \in \Plans(\Source,\Sink,\Domain)
  \right\}\; .
\end{equation}
The resulting infimum defines a distance $\Wass{1}(\Source,\Sink)$
between $\Source $ and $\Sink$, which is called the Wasserstein-1 or
Kantorovich-Rubinstein distance.

\subsection{The \MKEQS}
\label{sec:mkeqs} 
The above problem can be rewritten in many different forms (see
e.g.~\cite{Ambrosio:2003}). Here we are interested in the PDE
formulation, called \emph{Monge-Kantorovich equations} (\MKEQS),
introduced in~\cite{Evans-Gangbo:1999}. We present these equations as
formulated in~\cite{Bouchitte-Buttazzo:2001}.  Denoting with
$\Lip_{1}(\Domain)$ the space of continuous functions with Lipschitz
constant equal to $1$, we look for the pair
$(\OptPot,\OptTdens) \in (\Lip_{1}(\Domain),\RadonPlus(\Domain))$ that
satisfies the following nonlinear system of PDEs:
\begin{subequations}
  \label{eq:mkeqs}
  \begin{align}
    \label{eq:mkeqs-elliptic}
    -\Div[\Tdens] (\Tdens \Grad[\Tdens] \Pot)
    = 
    \Source-\Sink
    \quad 
    &\textrm{on $\Domain$}
    \\
    \label{eq:mkeqs-gradless1}
    \Norm{\Grad[\Tdens]\Pot} \leq 1 
    \quad &\textrm{on $\Domain$}
    \\
    \label{eq:mkeqs-eikonal}
    \Norm{\Grad[\Tdens]\Pot} = 1 \quad
    &\Tdens-a.e. \mbox{ in } \Domain\; ,
  \end{align}
\end{subequations}
where $\Grad[\Tdens]$ denotes the gradient with respect to the
measure $\Tdens$ as described in~\cite{Bouchitte-et-al:1996,
  Fragala-Mantegazza:1999}.  Note that \cref{eq:mkeqs-elliptic} must
be interpreted in the following weak form:
\begin{equation}
  \label{eq:mkeqs-weak}
  \int_{\Domain} 
  \Grad[\Tdens]  \Pot(x) \cdot \Grad[\Tdens] \Test(x)
  \mathrm{d}\Tdens (x)
  = 
  \int_{\Domain}
  \Test(x)
  \mathrm{d} \Source (x)
  - 
  \int_{\Domain}
  \Test(x)
  \mathrm{d} \Sink (x)
  \quad
  \forall \Test \in \Cont[1](\Domain)\; ,
\end{equation}
which implies zero Neumann boundary conditions and must be
complemented by the constraint that $\Pot$ has zero average.  The
components of the solution pair $(\OptPot,\OptTdens)$ of the above
system are named Kantorovich potential and \OTD. Among the several
results available on the \MKEQS, we recall here the few that will be
used in this work (see, e.g., ~\cite{Feldman-McCann:2002,
  Ambrosio:2003, DePascale-et-al:2004, Santambrogio:2009}).
Specifically, if either one of the measures $\Source$ and $\Sink$
admits a density with respect to the Lebesgue measure, then
$\OptTdens$ admits a unique density with respect to the Lebesgue
measure. Moreover, when $\OptTdens$ admits a density the notion of
gradient $\Grad[\Tdens]$ with respect to a measure coincides with the
notion of Sobolev weak gradient (see \cite[Example 2.3,
  Appendix]{Bouchitte-et-al:1996}). It is important to remark that
there exist infinitely many functions $\OptPot$ that
satisfy~\cref{eq:mkeqs}. In fact outside $\Supp(\OptTdens)$ it is
possible to perturb any solution $\OptPot$ and still satisfy all
constraints of~\cref{eq:mkeqs}.  However, in the support of
$\OptTdens$ the Kantorovich potential is
unique~\cite{Evans-Gangbo:1999}.

\subsection{The Beckmann Problem}
\label{sec:beckmann}
The \MKEQS\ are related also to another equivalent formulation of the
$\Lspace{1}$-\OTP, called \emph{Beckmann Problem}
(see~\cite{Santambrogio:2015}). The latter problem tries to find the
optimal vector-valued measure $\Vel\in[\Radon(\Domain)]^{\Dim}$ that
solves the following minimization problem:
\begin{equation}
  \label{eq:beckmann}
  \min_{\Vel \in [\Radon(\Domain)]^{\Dim}} 
  \int_{\Domain} \ABS{d\Vel} 
  \; : \;
  \Div(\Vel) = \Source-\Sink \; ,
\end{equation}
where $\ABS{d\Vel}$ denotes the total variation measure of $d\Vel$
(see~\cite{Santambrogio:2015}). The \MKEQS\ and the Beckmann Problem
are equivalent in the sense that $\OptVel$ is related to $\OptTdens$
and $\OptPot$ by:
\begin{equation}
  \label{eq:beckmann-mkeqs}
  d\OptVel = -d\OptTdens \Grad_{\OptTdens} \OptPot \;.
\end{equation}
Note that, when at least one of $\Source$ or $\Sink$ admits
$\Lspace{1}$-density with respect to the Lebesgue measure, then
$\OptVel=-\OptTdens\Grad\OptPot\in\Vof{\Lspace{1}(\Domain)}{\Dim}$ and
problem~\cref{eq:beckmann} can be rewritten as:
\begin{equation*}
  \min_{\Vel \in \Vof{\Lspace{1}(\Domain)}{\Dim}} 
  \int_{\Domain} |\Vel |\dx 
  \; : \;
  \Div(\Vel) = \Source-\Sink \; .
\end{equation*}

\subsection{The Dynamic Monge-Kantorovich equations}
\label{sec:DMK}

The model proposed in~\cite{Facca-et-al:2018,Facca-et-al-numeric:2020}
considers the case of either $\Source$ or $\Sink$ having
$\Lspace{2}$-density with respect to the Lebesgue measure. With an
abuse of notation, we will denote with $\Source$ and $\Sink$ both the
measures and their densities.  The model, named Dynamic
Monge-Kantorovich (\DMK) formulation (encapsulated here
  within a separate problem statement as it forms the basic starting
  point for our numerical developments), reads as follows.
\begin{Problem}[Dynamic Monge-Kantorovich (\DMK) equations]
  \label{prob:dmk}
  Consider a convex and bounded domain $\Domain$ in $\REAL^\Dim$ with
  smooth boundary and let $\Source,\Sink:\Domain\rightarrow\REALNN$ be
  such that $\int_{\Domain} \Source(x) \dx = \int_{\Domain} \Sink
  (x)\dx $. Find the pair $(\Tdens,\Pot)$ with $\Tdens:\REALNN\times
  \Domain \mapsto \REALNN$ and $\Pot:\REALNN\times \Domain \mapsto
  \REAL$ that satisfies:
  \begin{equation}
    \label{eq:flatmodel}
    \left\{
      \begin{aligned}
	&- \Div \big(\Tdens(t,x)\Grad \Pot(t,x)\big)=\Source(x)-\Sink(x)\\
	&\Dt{\Tdens}(t,x) = \Tdens (t,x) \big(|\Grad \Pot(t,x)|-1\big)\\
	&\Tdens(0,x) = \Tdens_0(x) >0.
      \end{aligned}
    \right.
  \end{equation}
  with zero Neumann boundary conditions on $\partial\Domain$ and $\int_{\Domain}\Pot=0$.
\end{Problem}
In \cite{Facca-et-al:2018} the authors conjecture that the solution
$(\Tdens(t,\cdot), \Pot(t,\cdot))$ of~\cref{prob:dmk} is
asymptotically equivalent to the solution $(\OptTdens, \OptPot)$ of
the \MKEQS~\cref{eq:mkeqs}. In support of this conjecture,
in~\cite{Facca-et-al-numeric:2020} the authors propose a
Lyapunov-candidate functional $\Lyap$ defined on non-negative
densities $\Gen{\Tdens}\in\Lspace{1}(\Domain)$ and given by the sum of
an energy and a mass functional, i.e.:
\begin{gather}
  \label{eq:lyap-def}
  \Lyap(\Gen{\Tdens})
  := \Ene(\Gen{\Tdens}) + \Wmass(\Gen{\Tdens})
  \\
  \label{eq:lyap-ene-wmass}
  \Ene(\Gen{\Tdens})
  := \sup_{\Test \in \Lip{}(\Domain)}
  \left\{
    \int_{\Domain} \left(
      (\Source-\Sink)\Test-\Gen{\Tdens}\frac{\ABS{\Grad\Test}^2}{2}
    \right)\dx
  \right\}
  \quad
  \Wmass(\Gen{\Tdens})
  :=\frac{1}{2}\int_{\Domain}{\Gen{\Tdens}\dx} \;.
\end{gather}
The functional $\Lyap$ can be rewritten as:
\begin{equation*}
  \Lyap(\Gen{\Tdens})=
  \frac{1}{2}\int_{\Domain}{\Gen{\Tdens}\ABS{\Grad \PotOp(\Gen{\Tdens})}^2}\dx
  + 
  \frac{1}{2}\int_{\Domain}{\Gen{\Tdens}\dx} 
\end{equation*}
as long as the energy $\Ene$ in~\cref{eq:lyap-ene-wmass} admits a
maximizer $\PotOp(\Gen{\Tdens} )$ that is the weak solution of the PDE
$-\Div(\Gen{\Tdens}\Grad\Gen{\Pot})=\Source-\Sink$ with zero-Neumann
boundary condition.

The main advantage of this formulation is that the introduction of the
dynamics allows the implementation of efficient numerical algorithms
for the solution of the original MK
equations~\cite{Facca-et-al-numeric:2020}.  We summarize in the
following Proposition the main results proved
in~\cite{Facca-et-al:2018, Facca-et-al-numeric:2020}.
\begin{Prop}
\label{prop:lyap}  
The \OTD\ $\OptTdens$ is the unique minimizer of $\Lyap$ and the
corresponding infimum value is exactly the Wasserstein-1 distance
between $\Source$ and $\Sink$:
$\Lyap(\OptTdens)=\Wass{1}(\Source,\Sink)$.  Moreover, for any $T>0$
for which a solution pair $(\Tdens(t,\cdot), \Pot(t,\cdot)$
of~\cref{prob:dmk} exists, for $t\in[0,T[$ the functional $\Lyap$
decreases along the $\Tdens$-trajectory.
\end{Prop}
More recently, in~\cite{Piazzon-et-al:2019} the authors addressed a
modified version of the \DMK\ problem and were able to prove global
existence of a time-asymptotic solution and its convergence (in a weak
sense) toward the solution of the \MKEQS.  The proof is based on the
reformulation of the dynamical system in~\cref{eq:flatmodel} as a
Gradient Flow in metric spaces~\cite{Ambrosio-et-al:2005} for the
functional $\Lyap$, contributing to strengthen the theoretical
background of the original \DMK\ model. On the other hand,
experimental numerical results with this modified formulation show
that it is not as efficient and robust as the original formulation
proposed in~\cite{Facca-et-al:2018, Facca-et-al-numeric:2020}.
For this reason, we continue working in this paper with the latter.

\subsection{Extension to manifolds}
\label{sec:otp-manifold}
We restrict our discussion to the case of a compact manifold
$(\Manifold,\Metric)$ ($\Manifold$ for short) with no boundary and
equipped with a smooth metric $\Metric$.  We denote with
$\Scal{v}{w}_{\Metric[\Point]}$ the application of the metric
$\Metric$ evaluated at $\Point$ to the two vectors
$v,w\in\Tangent[\Point]{\Manifold}$.  Moreover we denote with
$\Grad[\Metric]$ and $\Norm[\Metric]{\cdot}$ the gradient and the
vector norm with respect to the metric tensor $\Metric$. The distance
induced by the metric $\Metric$ is defined as:
\begin{equation*}
 \Dist_{\Metric}(x,y) = \inf
 \left\{
   \int_{0}^{1}
   \sqrt{
     \Scal{\dot{\Curve}(s)}{\dot{\Curve}(s)}_{\Metric[\Curve(s)]}
   } 
   \; \mathrm{d}s
   \;:\; 
   \begin{gathered}
     \Curve\in \Cont[1]([0,1],\Manifold)\\
     \Curve(0) = x, \Curve(1) =y
   \end{gathered}
 \right\},
\end{equation*}
and we denote with $\dv$ the volume form induced by the metric
$\Metric$.

Now we can proceed with the formulation of the \OTP\ on a manifold
$(\Manifold,\Metric)$, adapting the notation and results
of~\cite{Pratelli:2005} to our setting.
Replacing in~\cref{eq:otp-plan,eq:l1-otp} the domain $\Domain$ with
$\Manifold$ and the Euclidean distance term $|x-y|$ with the
distance $\Dist_{\Metric}(x,y)$, the $\Lspace{1}$-\OTP\ on $\Manifold$
tries to find the optimal plan
$\Opt{\Plan}\in \Plans(\Source,\Sink,\Manifold)$ that solves the
minimization problem:
\begin{equation*}
  \inf_{\Plan}
  \left\{
    \int_{\Manifold \times \Manifold} \Dist_{\Metric}(x,y) \mathrm{d}\Plan(x,y)
    \ : \ 
    \Plan \in \Plans(\Source,\Sink,\Manifold)
\right\} \; .
\end{equation*}
This infimum value is exactly the Wasserstein-1 distance 
$\Wass{1,g}(\Source,\Sink)$ on $\Manifold$.
The Monge-Kantorovich equations described in~\cref{eq:mkeqs} can be
generalized as well and become:
\begin{subequations}
  \label{eq:mkeqs-manifold}
  \begin{align}
    \label{eq:mkeqs-manifold-elliptic}
    -\Div[\Metric,\Tdens] (\Tdens \Grad[\Metric,\Tdens] \Pot)
    &= 
    \Source -\Sink 
    \quad 
    &\textrm{on $\Manifold$}
    \\
    \label{eq:mkeqs-manifold-gradless1}
    \Norm[\Metric]{\Grad[\Metric,\Tdens]\Pot} &\leq 1 
    \quad &\textrm{on $\Manifold$}
    \\
    \label{eq:mkeqs-manifold-eikonal}
    \Norm[\Metric]{\Grad[\Metric,\Tdens]\Pot} &= 1 \quad
    &\Tdens-a.e.
  \end{align}
\end{subequations}
where \cref{eq:mkeqs-manifold-elliptic} must be interpreted in 
the following weak form
\begin{multline}
  \label{eq:mkeqs-weak-manifold}
  \int_{\Manifold} 
  \Scal{
    \Grad[\Metric,\Tdens]\Pot(\Point)
  }{
    \Grad[\Metric,\Tdens]\Test(\Point)}_{\Metric(\Point)}
  \mathrm{d}\Tdens(\Point)
  =\\= 
  \int_{\Manifold}
  \Test(\Point)
  \mathrm{d}\Source(\Point)
  - 
  \int_{\Manifold}
  \Test(\Point)
  \mathrm{d}\Sink(\Point)
  \quad
  \forall \Test \in \Cont[1](\Manifold) \; ,
\end{multline}
where $\int_{\Manifold}\Pot=0$ is assumed. Likewise the Euclidean
case, if $\OptTdens$ admits a density with respect to the measure
$\dv$, the gradient $\Grad[\Metric,\Tdens]$
in~\cref{eq:mkeqs-manifold} coincides with the classical weak gradient
$\Grad[\Metric]$ on $\Manifold$.  Existence and uniqueness results of
the \OTD\ in a Riemannian setting were proved
in~\cite{Feldman-McCann:2001}, while $\Lspace{p}$ summability is still
an open question.

Beckmann Problem (\cref{eq:beckmann}) can be transported as well into a
manifold becoming:
\begin{equation}
  \label{eq:beckmann-manifold}
  \min_{\Vel \in [\Radon(\Manifold)]^{\Mdim}} 
  \int_{\Manifold} \Norm{d \Vel }_{\Metric} 
  \; : \;
  \Div_{\Metric}(\Vel) = \Source-\Sink
\end{equation}
and~\cref{eq:beckmann-mkeqs} holds also in this Riemannian setting.
Again, when at least one of $\Source$ or $\Sink$ admits
$\Lspace{1}$-density with respect to the volume form,
~\cref{eq:beckmann-manifold} can be rewritten as:
\begin{equation*}
  \min_{\Vel \in \Vof{\Lspace{1}(\Manifold)}{\Mdim}} 
  \int_{\Manifold} |\Vel |_{\Metric}\dv 
  \; : \;
  \Div_{\Metric}(\Vel) = \Source-\Sink \; ,
\end{equation*}
where, similarly to~\cref{eq:beckmann-mkeqs}, the tangent velocity
field is given by:
\begin{equation}
  \label{eq:beckmann-mkeqs-manifold}
  \OptVel = -\OptTdens\Grad[\Metric]\OptPot \;.
\end{equation}

Finally, we discuss next the extension of the \DMK\ formulation from
the Euclidean to the Riemannian setting.  This is obtained by simply
replacing the Euclidean differential operators with the Riemannian
ones identified with the metric $\Metric$ as subscript.  Thus we will
write $\Norm{\ }_{\Metric}$, $\Scal{}{}_{\Metric}$, $\Div_{\Metric}$,
and $\Grad[\Metric]$ for $\Norm{\ }$, $\cdot$, $\Div$, and $\Grad$,
respectively.  Then, the Riemannian version of~\cref{prob:dmk} can be
written as follows.
\begin{Problem}[Surface Continuous \DMK]
  \label{prob:dmk-surface} 
  Let $(\Manifold,\Metric)$ be a $\Mdim$-dimensional
  smooth compact manifold  with no
  boundary and let $\Source,\Sink\in\Lspace{2}(\Manifold)$ be such
  that $\int_{\Manifold}\Source\dv=\int_{\Manifold}\Sink \dv$. Find
  the pair  $(\Tdens,\Pot)$ with $\Tdens:\REALNN\times
    \Manifold\mapsto \REALNN$ and $\Pot:\REALNN\times \Manifold\mapsto
    \REAL$ satisfying:
  \begin{equation}
    \label{Eq:surfmodel}
    \left\{
      \begin{aligned}
        &
        - \Div[\Metric]
        \left(\Tdens(t,\Point)\Grad[\Metric]\Pot(t,\Point)\right)
        =\Source(x)-\Sink(x) \\
        &
        \Dt{\Tdens}(t,\Point)=\Tdens(t,\Point)
        \left(\Norm{\Grad[\Metric]\Pot(t,\Point)}-1\right)\\
        &\Tdens(0,\Point) = \Tdens_0(\Point)>0\; ,
      \end{aligned} \right.
  \end{equation}
  with $\Pot$ having zero mean.
\end{Problem}
The generalization to the Riemannian setting of the \LCF\
in~\cref{eq:lyap-def} for any non-negative
${\Tdens}\in \Lspace{1}(\Manifold)$ is also straightforward:
\begin{gather}
  \label{eq:lyap-def-manifold}
  \Lyap_{\Metric}({\Tdens})
  := \Ene_{\Metric}({\Tdens}) + \Wmass_{\Metric}({\Tdens})
  \\
  \nonumber
  \Ene_{\Metric}({\Tdens})
  := \sup_{\Test \in \Lip{}(\Manifold)}
  \left\{
    \int_{\Manifold} \left(
      (\Source-\Sink)\Test-{\Tdens}\frac{\ABS{\Grad[\Metric]\Test}^2}{2}
    \right)\dv
  \right\}
  \quad
  \Wmass_{\Metric}({\Tdens})
  :=\frac{1}{2}\int_{\Manifold}{{\Tdens}\dv} \;.
\end{gather}
It is immediate to show that, mutata mutandis, \Cref{prop:lyap} holds
also in this setting. It is thus natural to extend also the conjecture
of the convergence of the solution of~\cref{Eq:surfmodel} toward the
solution of the \MKEQS\ in~\cref{eq:mkeqs-manifold}.

\section{Numerical Discretization}
\label{sec:numerical-discretization}

This section is dedicated to the development of the numerical
discretization for \cref{prob:dmk-surface} obtained by extending to
the surface setting the numerical DMK formulation proposed
in~\cite{Facca-et-al:2018, Facca-et-al-numeric:2020} for the Euclidean
framework. The approach is based on the Method Of Lines (MOL) and
considers a 2-dimensional compact surface $\Surf$ with no boundary,
embedded in $\REAL^3$.
\nameCrefs{Eq:surfmodel}~\labelcref{Eq:surfmodel} are first
discretized in space and then the resulting ODE system is discretized
in time.  The long-time numerical solution of the latter is used
to approximate the solution of the \MKEQS.
We start by describing the (formal) weak formulation of the continuous
surface \DMK\ system obtained by exploiting the Hilbert space
structure of $\Hilb{1}(\Surf)$ and $\Lebesgue{2}(\Surf)$.  We test the
system using functions in these spaces to yield:
\begin{subequations}
  \label{eq:MKCont}
  \begin{align}
    &\int_{\Surf}\Tdens(t) \Scal{\Grad[\Metric]
      \Pot(t)}{\Grad[\Metric]\testp}_{\Metric}\dv
    =\int_{\Surf}(\Source-\Sink)\testp\dv, & \testp \in
    \Hilb{1}(\Surf) \label{eq:MKCont-1} \\
    &\int_{\Surf}\Dt{\Tdens}(t)\testx(x)\dv\!=\!
    \int_{\Surf}
      \Tdens(t)\left(\Norm{\Grad[\Metric]\Pot(t)}_{\Metric}-1
    \right)\testx
    \dv ,
    &  \testx \!\in \Lebesgue{2}(\Surf) \label{eq:MKCont-2}\\
    &\int_{\Surf}\Tdens(0)\testx\dv =
    \int_{\Surf}\Transport_0\testx\dv  & \testx \in
    \Lebesgue{2}(\Surf) \label{eq:MKCont-3} \\
    &\int_{\Surf}\Pot\testp\dv =0 & \testx \in
    \Hilb{1}(\Surf)\, , \label{eq:MKCont-4}
  \end{align}
\end{subequations}
ho aggiunto l'equazione per la media where, for simplicity, we have
dropped the spatial dependence of the dependent variables.

\subsection{Spatial discretization}
\label{sec:spatial-discretization}

Spatial discretization of~\cref{eq:flatmodel} follows the Surface
Finite Element Method (SFEM) described in~\cite{Dziuk-Elliott:2013}.
SFEM can be summarized in two steps: i) discretization of the surface
and ii) definition of the finite dimensional FEM spaces.
In the first step, we assume that the smooth surface $\Surf$ can be
subdivided by a surface triangulation $\Triang(\Surf)$ formed by the
union of non-intersecting surface geodesic elements (triangles)
$\Cell[i]\in\Surf$, $i=1,\ldots,\NCell$.  Thus, similarly to the flat
case, we have that $\Surf=\Triang(\Surf)=\cup_{i=1}^{\NCell}\Cell[i]$
and $\Edge[ij]=\Cell[i]\cap\Cell[j]$ is a geodesic curve on $\Surf$
connecting triangle vertices $i$ and $j$.
This triangulation is then approximated by its piecewise linear
interpolant $\TriangH(\Surf)=\Surfh$, where $\Surfh$ is
defined by the union of 2-simplices in $\REAL^3$ (flat
three-dimensional triangles) having the same vertices of
$\Triang(\Surf)$.
The triangulation $\TriangH(\Surf)$ is assumed to be closely inscribed
in the sense of~\cite{Morvan}, or equivalently, in the sense
of~\cite{Dziuk-Elliott:2013}. Thus, we require that
$\TriangH(\Surf)\subset\Neigh[\epsilon]$, where
$\Neigh[\epsilon]\supset\Surf$ is a tubular neighborhood of $\Surf$ of
radius $\epsilon$ such that every point $\Point\in\Neigh[\epsilon]$
has a unique orthogonal projection onto $\Surf$.
Given a triangle $\CellH\in\TriangH(\Surf)$, we can directly extend the
classical two-dimensional ``flat'' definitions characterizing proper
triangulations. Thus, we denote the mesh parameter and mesh radius
as $\MeshPar=\max_{\Cell\in\TriangH}\MeshPar[\Cell]$ and
$\Inradius=\min_{\Cell\in\TriangH}\Inradius[\Cell]$, respectively, where
$\MeshPar[\Cell]$ is the length of the longest edge of $\CellH$ and
$\Inradius[\Cell]$ the radius of the circle inscribed in $\CellH$.
Moreover, we assume that $\TriangH(\Surf)$ is shape-regular, i.e.,
there exists a strictly positive constant $\rho>0$ independent of
$\MeshPar$ such that:
\begin{equation*}
  \frac{\Inradius[\Cell]}{\MeshPar[\Cell]}\ge\rho
  \qquad \forall\CellH\in\TriangH(\Surf)\;.
\end{equation*}
\begin{Remark}
  As a consequence of the previous definitions, for any flat cell
  $\CellH\subset\Surf[\MeshPar]$ there is a unique curved cell
  $\Cell\subset\Surf$, and this correspondence is bijective.
  Hence, all the properties of $\TriangH(\Surf)$ are inherited
  by $\Triang(\Surf)$ and can be given indifferently for only
  one of the two triangulations.
\end{Remark}
%

Following the approach described in~\cite{Facca-et-al-numeric:2020},
the two dependent variables are defined in two different
triangulations and two different functional spaces to avoid
oscillations. This sort of inf-sup stability requirement is not yet
fully understood and was used in
both~\cite{Facca-et-al:2018,Facca-et-al-numeric:2020} as a remedy for
observed checkerboard-like oscillations. Thus, we employ the
triangulation $\Triang[\MeshPar/2](\Surf)$ generated from
$\Triang[\MeshPar](\Surf)$ by conformally refining each flat triangle.
The new nodes in the conformal refinement are not moved back to the
surface $\Surf$ to enhance the stabilization properties introduced by
the refined mesh.
The sub-triangles in the refined triangulation
$\Triang[\MeshPar/2]$ are denoted by $\SubCell_{\Itdsubcell}$ while
$\Cell_{\Itdcell}$ identifies triangles in
$\Triang[\MeshPar](\Surf)$. We use
$\SubCell_{\Itdsubcell}\subset\Cell_{\Itdcell}$ to denote those cells of
$\Triang[\MeshPar/2](\Surf)$ that belong to $\Cell_{\Itdcell}$.

Given a smooth function $\Fgen:\Surf\to\REAL$, its smooth extension in
$\Neigh[\epsilon]$ is well defined and is denoted by
$\bar{\Fgen}:\Neigh[\epsilon]\to\REAL$.  This smooth extension is used
to define the tangential gradient of $\Fgen$ as
follows:
\begin{Defin}[Tangential Gradient]
  \label{def:emb-tan-gradient}
  Let $\bar{\Fgen}:\Neigh[\epsilon]\to\REAL$ be the continuous extension
  of a smooth function $\Fgen:\Surf\to\REAL$.
  The \emph{tangential gradient} of $\Fgen$ at a point
  $\Point\in\Surf$ is given by
  \begin{equation}\label{eq:gradgamma}
    \Grad[\Surf] \Fgen(\Point) = \Grad\bar{\Fgen}(\Point) -
    \Scal{\Grad\bar{\Fgen}(\Point)}{\Normal(\Point)}\Normal(\Point) = 
    \ProjMat(\Point)\Grad\bar{\Fgen}(\Point),
  \end{equation}
  is the unit oriented normal to the surface $\Surf$ at $\Point$, the
  projection tensor is given by
  $\ProjMat(\Point)=\IDtens-\Normal(\Point)\otimes\Normal(\Point)$,
  the bilinear form $\Scal{\cdot}{\cdot}$ is the standard scalar
  product in $\REAL^3$.
\end{Defin}
It turns out~\cite{Dziuk-Elliott:2013} that the tangential gradient is
equivalent to the metric gradient
(i.e. $\Grad[\Metric]\Fgen=\Grad[\Surf]\Fgen$) for any regular surface
$\Surf$. Moreover, Green's lemma holds for these tangential operators.

Then, we can define the SFEM DMK on $\Surfh$ by simply employing the
quantities defined in~\cref{def:emb-tan-gradient} and
projecting~\eqref{eq:MKCont} onto the finite dimensional
FEM spaces $\Vspace\subset\Hilb{1}(\Surfh)$ and
$\Wspace\subset\Lebesgue{2}(\Surfh)$.  Note that
$\NormalH=\Normal(\Point)|_{\Cell}$ is constant on each triangle
$\Cell\in\TriangH(\Surfh)$, and so is
$\ProjMat[\MeshPar]=\IDMat-\NormalH\otimes\NormalH$. Thus it is
convenient to use as test functions the Lagrange-polynomial basis
functions of $\Vspace$ and $\Wspace$, which are defined also
triangle-wise. Then, the integrals over $\Surfh$ transform into a sum
over $\TriangH(\Surfh)=\Surfh$, and we can write the following
SFEM-DMK semi-discrete problem embedded in $\REAL^3$:
\begin{Problem}
  Find $(\TdensH(t),\PotH(t))\in\Wspace\times\Vspace$ such that:
\begin{subequations}
  \label{eq:MKContH}
  \begin{align}
    &\int_{\Surfh}\TdensH(t) \Scal{\Grad[\Surf]
      \PotH(t)}{\Grad[\Surf]\testp_i}\dx
      =\int_{\Surfh}(\Source-\Sink)\testp_i\dx \; ,
    && i=1\ldots\Vdim \; ,
    \label{eq:MKCont-1H} \\
    &\int_{\Surfh}\Dt{\TdensH}(t)\testx_s\dx\!=\!
    \int_{\Surfh}
      \TdensH(t)\left(\Norm{\Grad[\Surf]\PotH(t)}-1
    \right)\testx_s    \dx \; ,
    &&  s=1,\ldots,\Wdim \; , \label{eq:MKCont-2H}\\
    &\int_{\Surfh}\Tdens(0)\testx_s\dx =
    \int_{\Surfh}\Transport_0\testx_s\dx \; ,
     && s=1,\ldots,\Wdim \; , \label{eq:MKCont-3H}\\
    &\int_{\Surfh}\PotH\testp_i\dx =0
     && i=1,\ldots,\Vdim \; , \label{eq:MKCont-4H}
  \end{align}
  where $\Vspace=\Span\{\testp_1,\ldots,\testp_{\Vdim}\}$ and
  $\Wspace=\Span\{\testx_1,\ldots,\testx_{\Wdim}\}$.
\end{subequations}
\end{Problem}
Note that we have transformed the surface integrals into a sum of
integrals on the planar triangles $\Cell$.  Similarly, we can define
the SFEM basis functions by means of standard (triangle-wise)
three-dimensional Lagrange interpolating polynomials and project their
Euclidean gradients onto each $\Cell\in\TriangH(\Surfh)$ using
$\ProjMat[\MeshPar]$.

Now we choose the same FEM spaces that worked in the planar
case~\cite{Facca-et-al-numeric:2020}. We denote by $\PONE(\SubCell)$
the set of affine functions supported on $\SubCell$ and with
$\PZERO(\Cell)$ the set of constant functions supported on $\Cell$.
Then we define $\Vspace=\PONE(\Triang[\MeshPar/2](\Surfh))$ and
$\Wspace=\PZERO(\Triang[\MeshPar](\Surfh))$, where
\begin{align*}
  &\PONE(\Triang[\MeshPar/2](\Surfh))=
  \left\{v\in\Cont[0](\Surfh)
      : v|_{\SubCell}\in\PONE(\SubCell), 
      \forall\SubCell\in\Triang[\MeshPar/2](\Surfh) 
    \right\} \; ,\\
  &\PZERO(\Triang[\MeshPar](\Surfh))=
  \left\{w
      : w|_{\Cell}\in\PZERO(\Cell),
      \forall\Cell\in\Triang[\MeshPar](\Surfh) 
    \right\} \; .
\end{align*}
Using FEM terminology, we choose to approximate the transport density
with non-conforming piecewise constant functions and the transport
potential with conforming piecewise affine functions.  Note that this
choice is compatible with the dynamic (second)
equation~\eqref{eq:MKCont-2H} since
$\Grad[\Surf]\PotH|_{\Cell}\in\PZERO(\Cell)$ if
$\PotH|_{\Cell}\in\PONE(\Cell)$.  Now, separating the temporal and
spatial variables, we can write the discrete transport density
$\TdensH(t,x)$ and the discrete Kantorovich potential $\PotH(t,x)$ and
its tangential gradient as:
\begin{gather*}
  \TdensH(t,x) =
  \sum_{\Itd=1}^{\muFEMdimension}
  \Tdens_\Itd(t)\testx_\Itd(x)\; ,
  \qquad
  \PotH(t,x) =
  \sum_{\Ipot=1}^{\uFEMdimension}
  \Pot_\Ipot(t)\testp_\Ipot(x)\; ,
  \\
  \Grad[\Surf]\PotH(t,x)=
  \sum_{\Ipot=1}^{\uFEMdimension} \Pot_\Ipot(t)\Grad[\Surf]
  \testp_\Ipot(x)\; .
\end{gather*}
We stress here that, by the way $\SubCell\in\Trianghh(\Surfh)$ is
built from $\Cell\in\TriangH(\Surfh)$,$\Cell$ is the union of four
subcells $\SubCell_{\ell}$, i.e., using a local enumeration,
$\Cell=\cup_{\ell=1}^{4}\SubCell_{\ell}$.  Thus, the computation of
the integral on the right-hand side of~\cref{eq:MKCont-2H} on
$\Cell\in\TriangH(\Surfh)$ requires a simple arithmetic average of the
gradients of the potential on the sub-triangles
$\SubCell\in\Trianghh(\Surfh)$:
\begin{equation*}
\int_{\Surfh}\TdensH\Norm{\Grad[\Surf]\PotH}\testx_{s}\dx=
\int_{\Cell[s]}\TdensH\Norm{\Grad[\Surf]\PotH}\dx=
\sum_{\ell=1}^{4}\int_{\SubCell[\ell]}\TdensH\Norm{\Grad[\Surf]\PotH}\dx\; .
\end{equation*}

Denoting with
$\Vect{\Tdens}(t)=\left\{\Tdens_\Itdcell(t)\right\}$
and
$\Vect{\Pot}(t)=\left\{\Pot_\Ipot(t)\right)\}$
the curves describing the time evolution of the spatially projected
finite-dimensional system, \cref{eq:MKCont} becomes:
\begin{subequations}
  \label{eq:dae}
  \begin{align}
    \label{eq:dae-ell}
    & \Matr{A}[\Vect{\Tdens}(t)]\ \Vect{\Pot}(t)=\Vect{b} \;,
      \\
    \label{eq:dae-ode}
    & \Dt{\Vect{\Tdens}}(t)=
      \Matr{D}(\Vect{\Pot}(t))\ \Vect{\Tdens}(t) \;,
      \qquad \Vect{\Tdens}(0)=\Vect{\Tdens_0} \; ,
  \end{align}
\end{subequations}
where the sparse stiffness matrix $\Matr{A}\Of{\Vect{\Tdens}(t)}$, and
the diagonal matrix $\Matr{D}\Of{\Vect{\Pot}(t)}$ are given by:
\begin{equation*}
  \Matr[i,j]{A}\Of{\Vect{\Tdens}(t)}=\int_{\Surfh} \TdensH (t)
  \Scal{\Grad[\Surf]\testp_i}{\Grad[\Surf]\testp_j}\dx
  \quad
  \Matr[\Itdcell,\Itdcell]{D}\Of{\Vect{\Pot}(t)}=
  \frac{1}{\ABS{{\Cell_{\Itdcell}}}}
  \int_{\Cell_{\Itdcell}}
  \left(
    \Norm{\Grad[\Surf]\PotH(t,x)}-1
  \right)
  \dx \, .
\end{equation*}
The zero mean constraint on $\PotH$ is incorporated in the linear
solution phase as detailed below.

\subsection{Time discretization}
\label{sec:time-discretization}

The numerical solution of \eqref{eq:dae} is obtained  by a forward Euler time
discretization. This allows an easy resolution of the nonlinearities
and an immediate decoupling of the resulting two algebraic systems of
equations.  Given a time sequence $(t_\tstep)$, ${k= 0, \ldots,\Kmax} $
with $t_{\tstepp}=t_{\tstep}+\Deltat[\tstep]$ we generate an
approximating sequence $(\TdensH^{\tstep},\PotH^{\tstep})$ defined as:
\begin{subequations}
  \nonumber
  \begin{align}
    & \Matr{A}[\Vect{\Tdens}^{\tstep}]\ \Vect{\Pot}^{\tstep}=\Vect{b} \;,
    \\
    \nonumber
    & \Vect{\Tdens}^{\tstepp} = \Vect{\Tdens}^{\tstep} 
    + \Deltat[\tstep] 
    \Matr{D}(\Vect{\Pot}^{\tstep})
    \Vect{\Tdens}^{\tstep}
  \end{align}
\end{subequations}
The first step involves the solution of a linear system whose matrix
$\Matr{A}[\Vect{\Tdens}^{\tstep}]$ is sparse and symmetric positive
semidefinite. This is obtained by means of the preconditioned
conjugate gradient algorithm using an IC(0)-based spectral
preconditioner with selective updates, as proposed
in~\cite{Facca-et-al-numeric:2020,Bergamaschi-et-al:2018} and a
deflation-like approach to incorporate the zero mean constraint.  The
second step involves the calculation of the cell gradients and a
direct time update that exploits the diagonal structure of the system.
At each time step, the step size $\Deltat[\tstep]$ is calculated on
the basis of the eigenvalues of matrix
$\Matr{D}(\Vect{\Pot}^{\tstep})$ to ensure stability of the forward
Euler time-stepping.

Time-convergence is considered achieved by the discrete solution
$(\TdensH^{\tstepp},\PotH^{\tstepp})$ when the relative variation of
$\TdensH$ is smaller than a given tolerance $\TolTime$, i.e.,
\begin{equation*}
  \Var(\TdensH):=
  \frac{
    \|\TdensH^{\tstepp}-\TdensH^{\tstep}\|_{\Lspace{\Perr}(\Surf)}
  }{
   \Deltat \|\TdensH^{\tstep}\|_{\Lspace{\Perr}(\Surf)}
 }<\TolTime
 \, .
\end{equation*}
We indicate with $\OptTime$ the time when equilibrium is numerically
reached and with $\OptTdensH$ and $\OptPotH$ the corresponding
$\TdensH^{\tstep}$ and $\PotH^{\tstep}$.

\subsection{Approximate solutions of the MK equations and
  Beckmann problem}
\label{sec:approximate-solutions}

The time-converged solutions $\OptTdensH$ and $\OptPotH$ are the
approximations of the solutions to the MK equations $\OptTdens$ and
$\OptPot$.  From these, we can calculate the approximate solution to
Beckmann Problem given in~\cref{eq:beckmann-mkeqs-manifold}.  Note
that the gradient of $\PotH(t)$ is defined on the refined triangles of
$\Triang[\MeshPar/2](\Surfh)$ while $\TdensH(t)$ is defined on the
triangles of $\Triang[\MeshPar](\Surfh)$.  As done in the numerical
solution algorithm, we project the gradient of $\PotH(t)$ onto
$\Triang[\MeshPar](\Surfh)$ by simple averaging. Our approximate
velocity $\VelH(t)\in\Vof{\PZERO(\Triang[\MeshPar](\Manifold))}{3}$ is
then given by:
\begin{equation} 
  \label{eq:velh}
  \VelH(t) |_{\Cell[\Itdcell]} := 
  \TdensH(t) |_{\Cell[\Itdcell] }
  \frac{1}{4}
  \sum_{\SubCell[\Itdsubcell] \in \Cell[\Itdcell] }
  \Grad[\Surf] \PotH(t) |_{\SubCell[\Itdsubcell]}
\end{equation}
We will denote $\OptVelH$ the velocity at the equilibrium time
$\OptTime$.

In addition, we can calculate the Wasserstein-1 distance
$\Wass{1}(\Source,\Sink)$ between $\Source$ and $\Sink$.
According to~\cref{prop:lyap}, $\Wass{1}(\Source,\Sink)$
corresponds to the minimum of the functional
$\Lyap_{\Metric}$, which can be estimated as:
\begin{equation*}
  \Wass{1}^{\MeshPar}(\Source,\Sink):= \Lyap_{\Metric,\MeshPar}(\OptTdensH)=
    \frac{1}{2}\int_{\Surfh}\OptTdensH\ABS{\Grad[\Surf]\OptPotH}^2\dx
  + 
  \frac{1}{2}\int_{\Surfh}\OptTdensH\dx \;.
\end{equation*}

\section{Numerical experiments}
The test case devised for the verification of the temporal and spatial
convergence properties of the scheme is settled on the unit sphere
$\Sphere{2}$. A constant density supported on a geodesic quadrilateral
is rigidly moved from one position to another.  In this case, the
explicit formulas of the solution of the \MKEQS\ can be found. This
allows us to test the accuracy of the proposed numerical scheme by
computing the relative error in the approximation of the \OTD, the
optimal velocity field $\OptVel$ solution
of~\cref{eq:beckmann-manifold}, and the Wasserstein-1 distance.  We
also compare our numerical approximations with the ones obtained using
the approach described in~\cite{Solomon-et-al:2014}. In this case, a
quantitative comparison is done only in terms of accuracy, while we
only report the CPU time required by our method because of the
drastically different implementations of the companion software
packages.

\subsection{Comparison with closed-form solutions}
\label{sec:exact-optdens}

\begin{figure}
  \centerline{ 
    \includegraphics[width=0.45\textwidth,
    trim={0.0cm 0.55cm 1.7cm 0cm},clip
    ]{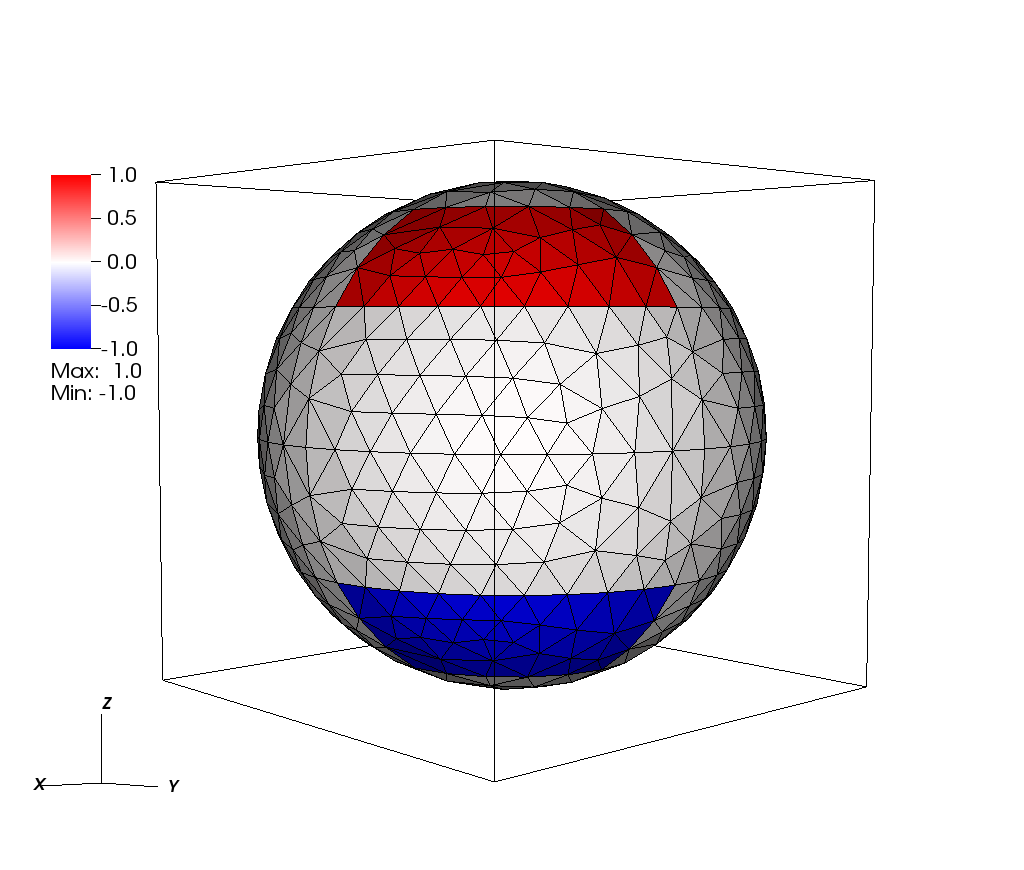}
    \quad
    \includegraphics[width=0.45\textwidth,
    trim={0.0cm 0.55cm 1.7cm 0cm},clip
    ]{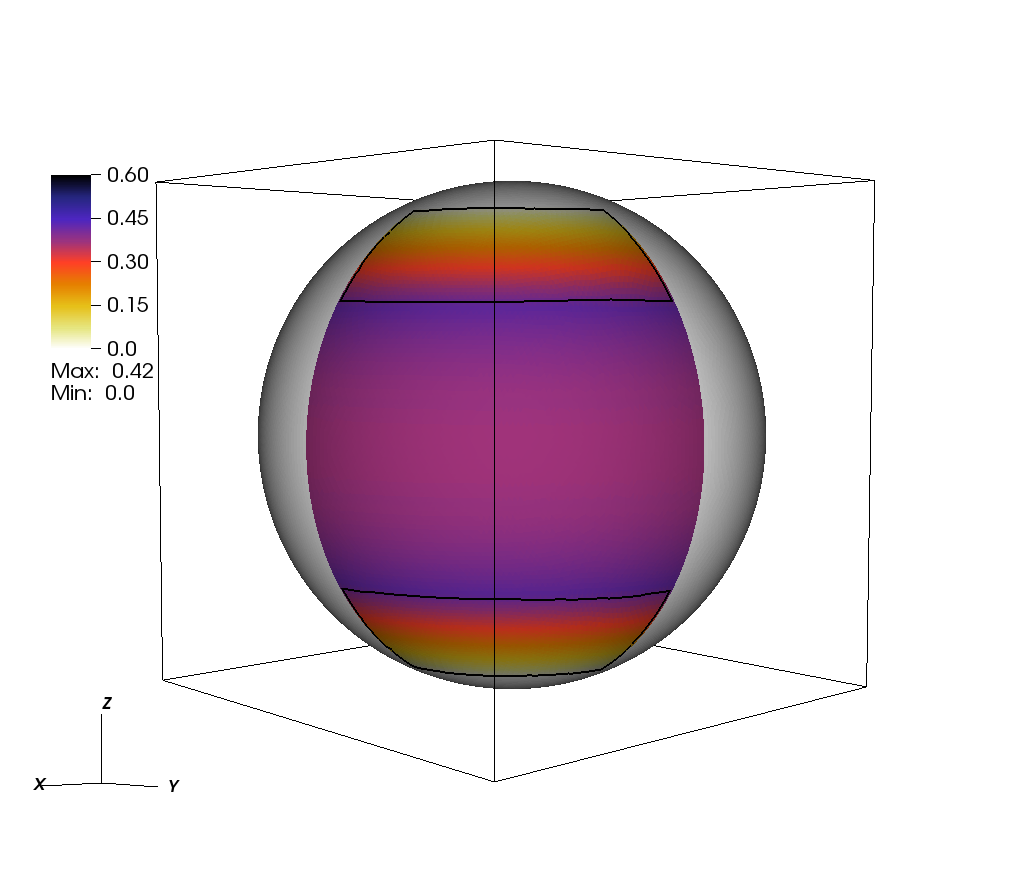}    
  }
  \caption{ Spatial distribution of the $\Lebesgue{2}$-projection of
    $\Source -\Sink$ on the mesh $\Triang[\MeshPar]$ with $554$ nodes
    and $1104$ triangles (left).  Spatial distribution of the density
    $\OptTdens$ together with the contour lines of the supports of
    $\Source$ and $\Sink$ (right).  }
  \label{fig:opt}
\end{figure}

\subsubsection{Description of the test case}

The test case considers the unit sphere $\Sphere{2}$ equipped with the
induced metric. The description of this test case is more easily
carried out using polar coordinates $(\SphereDist,\SphereAngle)$
defined with respect to the north pole $N=(0,0,1)$, where
$\SphereDist(\Point)=\Dist_{\Metric}(\Point,N)\in]0,\pi[$ and
$\SphereAngle(\Point)\in]0,2\pi[$ is the ``longitude'' angle.
We consider two unit densities $\Source$ and $\Sink$ with supports
given by:
\begin{align*}
  \Supp(\Source)&=\left\{(\SphereDist,\SphereAngle ) \ : \
    [\pi/6<\SphereDist<\pi/3] \quad \  \SphereAngle \in [0,\pi/2]\right\}
  \\
  \Supp(\Sink)&=\left\{(\SphereDist,\SphereAngle) \ : \
    [2\pi/3<\SphereDist<5\pi/6] \ \SphereAngle \in [0,\pi/2]\right\}\; .
\end{align*}
The approximate spatial distribution of $\Source$ and $\Sink$ is shown
in~\cref{fig:opt}.  For such densities, the pair $(\OptPot,\OptTdens)$
solution of the \MKEQS~\cref{eq:mkeqs-manifold} is given by:
\begin{gather*}
  \Opt{\tilde{\Tdens}}(\SphereDist,\SphereAngle)=
  \left\{
    \begin{aligned}
      & \left( cos(\pi/6) -\cos(\SphereDist)\right)/\sin(\SphereDist)
      && \mbox{if }\pi/6<\SphereDist<\pi/3 \  &&  0<\SphereAngle<\pi/2 
      \\
      & \left(cos(\pi/6) - cos(\pi/3)\right)/\sin(\SphereDist)
      && \mbox{if }\pi/3<\SphereDist<2\pi/3 && 0<\SphereAngle<\pi/2 
      \\
      & \left(cos(\pi/6) + cos(\SphereDist)\right)/\sin(\SphereDist)
      && \mbox{if } 2\pi/3<\SphereDist<5\pi/6 && 0<\SphereAngle<\pi/2
      \\
      & 0
      && \mbox{elsewhere}\\
    \end{aligned}
  \right.
  \\
  \Opt{\tilde{\Pot}}(\SphereDist,\SphereAngle)=-\SphereDist
\end{gather*}
We obtained these explicit formulas adapting to $\Sphere{2}$ the exact
solutions developed by~\cite{Buttazzo-Stepanov:2003} for the Euclidean
\MKEQS. We recall that $\Opt{\tilde{\Pot}}$ represents only one
possible solution of the \MKEQS, since in general the Kantorovich
potential is unique up to constants within the support of $\OptTdens$,
while outside it is not unique.  The spatial distribution of
$\OptTdens$ is shown on the right panel of~\cref{fig:opt}.
The explicit solution of Beckmann Problem can be derived
using~\cref{eq:beckmann-mkeqs}, obtaining:
\begin{equation}
 \label{eq:optvel}
  \Opt{\tilde{\Vel}}(\SphereDist,\SphereAngle)=\Opt{\tilde{\Tdens}}
  (\SphereDist,\SphereAngle)
  \left(
    1 \cdot \frac{\partial}{\partial \SphereDist},
    0 \cdot \frac{\partial}{\partial \SphereAngle }
  \right)\; ,
\end{equation}
where $\left(\frac{\partial}{\partial\SphereDist},
  \frac{\partial}{\partial\SphereAngle}\right)$ are the coordinate
vectors in $\Tangent[\CutPoint]{\Sphere{2}}$ with respect to the
coordinate system $(\SphereDist,\SphereAngle)$.

We measure the accuracy of our numerical approximation by calculating
the error in the solution of Beckmann problem $\Opt{\tilde{\Vel}}$ and
the error in the evaluation of the Wasserstein-1 distance
$\Wass{1}(\Source,\Sink)$.  The former is evaluated as:
\begin{equation*}
  \Err_{\BP}(\OptVelH):=
  \frac{\int_{\Surfh}|\OptVelH-\OptVel|\dv}
  {\int_{\Sphere{2}}|\OptVel|\dv}\; .
\end{equation*}
where the integral in the numerator is approximated via the mid-point
rule on $\Surfh$.
The use of the above error definition allows the evaluation of the
accuracy on both the velocity direction
$\Grad[\Surf]\Opt{\tilde{\Pot}}$ and magnitude given by
$\Opt{\tilde{\Tdens}}$.  Moreover, it is consistent with the
application of the approach proposed in
\cite{Solomon-et-al:2014} and thus it allows a fair comparison with
published results.

The error in the Wasserstein-1 distance between $\Source$ and $\Sink$
deserves a specific evaluation because of its importance in
applications. We then define:
\begin{equation*}
  \ErrWass(\Wass{1}^{\MeshPar}):=
  \frac{\ABS{\Wass{1}^{\MeshPar}(\Source,\Sink)-\Wass{1}(\Source,\Sink)}}
  {\Wass{1}(\Source,\Sink)} \; ,
\end{equation*} 
where $\Wass{1}^{\MeshPar}$ denotes the specific numerical
approximation of the Wasserstein-1 distance. For the test case
considered, the exact value can be calculated as:
\begin{equation*}
  \Wass{1}(\Source,\Sink) 
  = 
  \int_{\Sphere{2}}\OptTdens
  = 
  \left(
    \frac{\pi\sqrt{3}}{6}
    -\left(\sqrt{3}-1\right)\left(2-\frac{\pi}{3}\right)
  \right)
  \frac{\pi}{4}\approx
  0.876739625901484 \ .
\end{equation*}

Finally, we report a brief description of the implementation of the
\EMDADMM\ approach proposed by~\cite{Solomon-et-al:2014} , which is
our benchmark against which we compare our method.  In \EMDADMM, the
approximate solution $\ADMMOptVelH$ is obtained as the following
linear combination:
\begin{equation*}
  \ADMMOptVelH(x)=
  \Grad_{\Metric} \bar{\Pot}(x) +
  \sum_{i}^{\Neig} \Vect[i]{c} V_i(x) 
\end{equation*}
where $\bar{\Pot}$ solves the Poisson equation $-\Delta_{\Metric}
\bar{\Pot}=\Source-\Sink$, while $V_{i}(x)$ are vector fields built
from the piecewise linear approximation of the the first $\Neig$
eigenfunctions of the Laplace Beltrami operator (i.e., eigenvector of
the corresponding FEM stiffness matrix).  The optimal coefficient
vector $\Vect{c}$ is found by means of the Alternating Direction
Method of Multipliers (ADMM)~\citep{Boyd:2011}. In our tests we used
the suggested value $\Neig=N_{\MeshPar}/16$, where $N_{\MeshPar}$ is
the number of nodes in $\TriangH$.

\begin{figure}
  \centerline{
    \includegraphics[width=0.8\textwidth
    ]{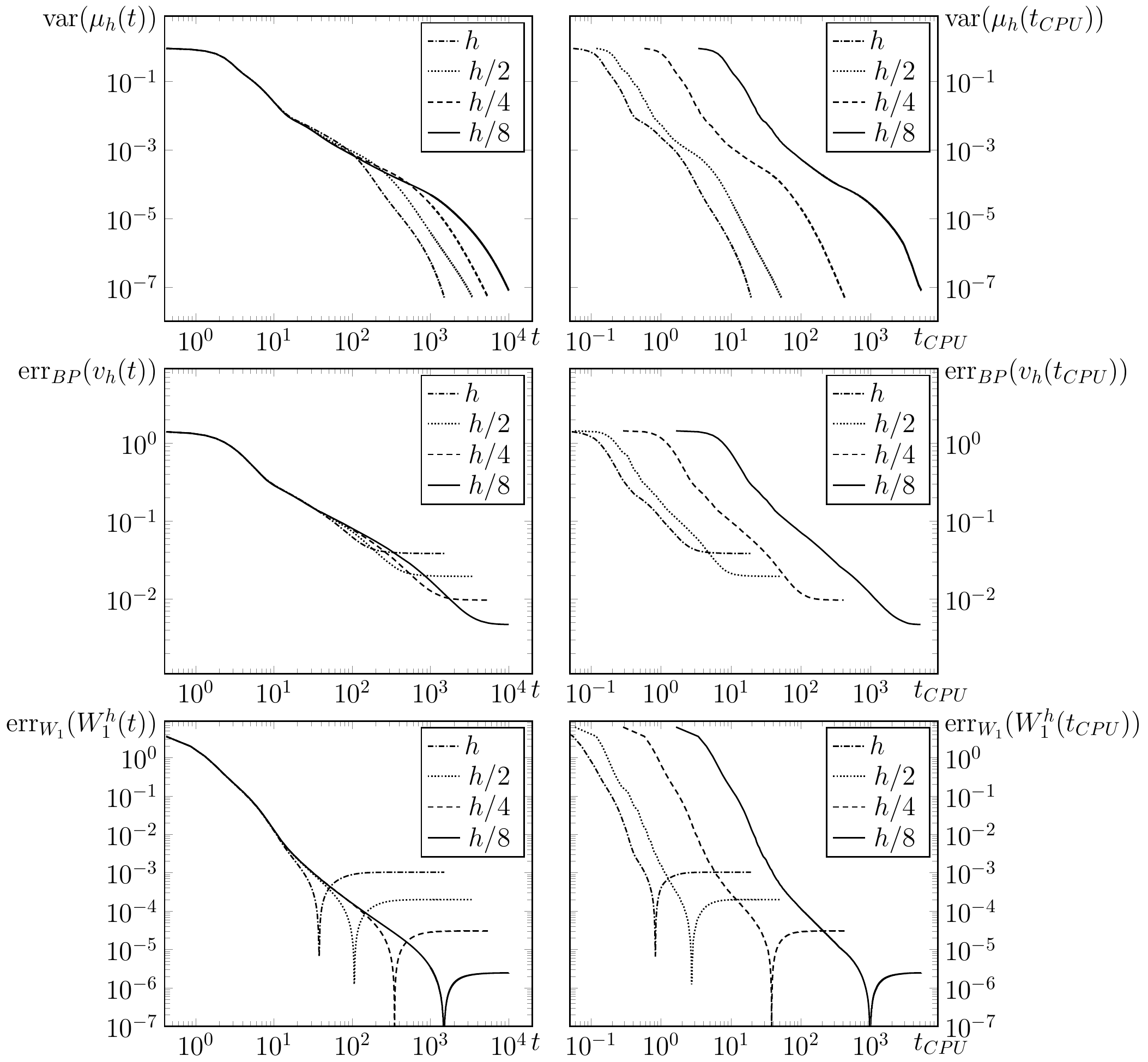}
  }
  \caption{ Experimental convergence towards the \MK\
    solution.  Left column: time ($t$) evolution of
    $\Var(\TdensH(t))$, which measures convergence towards
    equilibrium, and of $\Err(\VelH(t))$ and
    $\ErrWass(\TdensH(t))$, which measure the accuracy of
    the spatial approximation.  The results are obtained on
    successive uniform refinements of a an initial mesh with
    554 nodes and 1124 triangles.  Right column: values of
    the four simulation metrics as a function of
    computational time ($t_{CPU}$, seconds).  The
    simulations were conducted on a first-generation 2.2GHz
    Intel-I5 (1-core) laptop computer.}
  \label{fig:sphere}
\end{figure} 

In our experiments, $\Sphere{2}$ is approximated by an initial
triangulation ($\Sphere{2}_{\MeshPar}=\Surfh=\TriangH(\Sphere{2})$)
made up of 554 nodes and 1124 triangles (reported in~\cref{fig:opt})
obtained using the software described
in~\cite{Persson-Strang:2004}. This initial triangulation is
constrained to have nodes lying exactly on the meridians and parallels
defining the boundaries of $\Supp(\Source)$ and $\Supp(\Source)$.
Successive uniform refinements of the initial mesh $\Surfh$ are
obtained by calculating edge mid-points and moving them to the unit
sphere surface, taking care of preserving the alignment with
$\Supp(\Source)$ and $\Supp(\Sink)$. This process is repeated three
times, for a total of four mesh levels. At each refinement level, the
sub-grid $\Trianghh(\Sphere{2})$ is obtained by uniform refinement with
edge midpoints maintained in the original planar triangle and
not moved to the surface.

\subsubsection{Experimental results}

\paragraph{Convergence of Surface DMK.}
The experimental convergence results for the surface DMK are
summarized in~\cref{fig:sphere}.  Denoting with $\TdensH(t)$ the
piecewise linear interpolation of the sequence $\TdensH^\tstep$, each
row of~\cref{fig:sphere} show the time evolution of the quantities
$\Var(\TdensH(t))$, $\Err_{\BP}(\OptVelH(t))$, and
$\ErrWass(\TdensH(t))$.  Each plot contains the results obtained for
the four considered refinement levels. To better envision the
computational effort of our approach, the right panels report the same
quantities with respect to the CPU time $t_{CPU}(t)$, i.e., the CPU
time in seconds to arrive at simulation time $t$ on a first-generation
2.2GHz Intel-I5 (1-core) laptop computer.

The top-left panel in~\cref{fig:sphere} shows a monotone convergence
toward equilibrium for all refinement levels. The same monotonically
decreasing profile is observed in the $\Err_{\BP}(\VelH(t))$ plots.
Here the error saturates at a level that decreases by a constant
factor at every refinement, indicating the achievement of the spatial
resolution limit attainable with the given triangulation. In addition,
these results suggest first order convergence of the spatial
discretization for the solution of Beckmann problem.  Finally, the
bottom panels show the evolution of the relative error on the
Wasserstein-1 distance.  The non-monotone profiles can be explained as
follows. First of all, note that the Wasserstein-1 distance is
evaluated as the minimum of $\Lyap_{\Metric}$, and thus we are approximating
$\Wass{1}$ from above. On the other hand, the integrals are
approximated by the mid-point quadrature rule from below (the area of
the linear triangles is always smaller or equal than that of the
corresponding surface triangle). At a certain time $\hat{t}$, the
overestimation of $\Lyap_{\Metric}$ and the underestimation of the
integrals compensate and $\ErrWass(\TdensH(\hat{t}))$ becomes
artificially close to zero.  After $\hat{t}$, the error in the
approximation of the integral dominates and the global error saturates
at a value that depends on the mesh size $\MeshPar$.

\begin{figure}
  \centerline{
    \includegraphics[height=0.35\textwidth
    ]{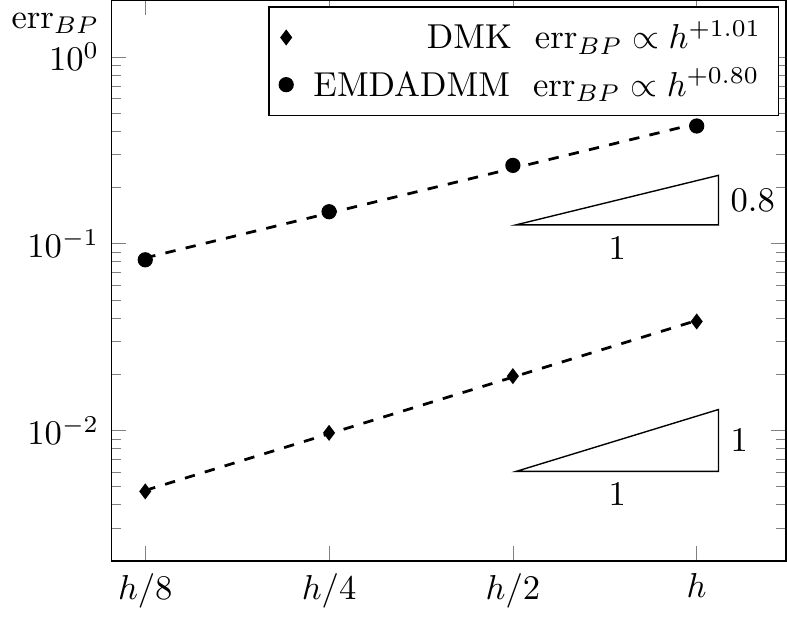}
    \quad
    \includegraphics[height=0.35\textwidth
    ]{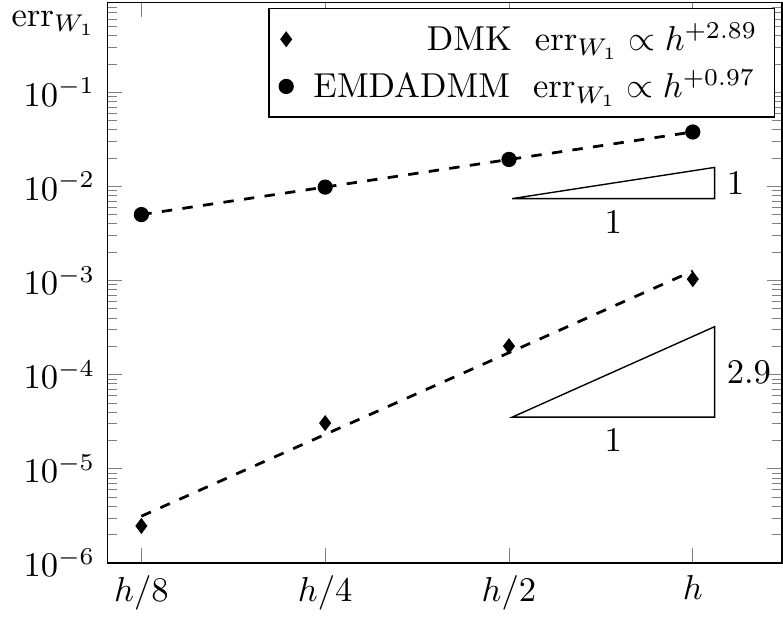}
  }
  \caption{Spatial convergence of \DMK\ and \EMDADMM\ on the relative
    Beckmann error $\Err_{\BP}(\OptVelH)$ (left) and on the
    Wasserstein-1 distance $\Err_{\Wass{1}}(\OptTdensH)$ (right).  The
    errors are calculated on four mesh refinements and the
    corresponding experimental spatial convergence rate with respect
    to mesh parameter $\MeshPar$, evaluated by linear approximation,
    is reported in the plot legend.}
  \label{fig:convergence}
\end{figure}

~\Cref{fig:convergence} shows the experimental mesh convergence
attained by the surface \DMK\ and the \EMDADMM\ methods on the
solution of the $\Lspace{1}$-OTP in terms of the velocity field of
Beckmann problem ($\Err_{\BP}(\OptVelH)$, left panel) and the
calculation of the Wasserstein-1 distance
($\Err_{\Wass{1}}(\OptTdensH)$, right panel).  Four nested mesh
refinements are used and the convergence lines (in log-log plot) are
calculated by linear least squares. In the first case the \DMK\
achieves a rate of order one, coherent with the piecewise constant
($\PONE$) discretization of the transport density and of the gradient
of the Kantorovich potential. This result is also consistent with the
results presented for the Euclidean case
in~\cite{Facca-et-al-numeric:2020}.  The experiments with the
\EMDADMM\ method show a slight sub-optimal rate of convergence of
about 0.8, providing a further indication of the lower accuracy of
this scheme with respect to \DMK.  The right panel
of~\cref{fig:convergence} shows the convergence plots when the two
schemes are applied for the calculation of the $\Wass{1}$ distance.
The experimental error obtained with \DMK\ scales almost
cubically with the mesh parameter $\MeshPar$, while \EMDADMM\
maintains first order convergence.

The accuracy attained by the DMK approach
is much higher than the one attained by the \EMDADMM\ method. Already
at the coarsest triangulation the DMK accuracy is well below $0.1\%$
with a CPU time of approximately 2 seconds.
In contrast, \EMDADMM\,
with a comparable computational effort, achieves on the the same mesh
an accuracy of about 4\%.

\section{Conclusion}
We propose the numerical solution of the $\Lspace{1}$-\OTP\ on
triangulated surfaces via the \DMK\ approach, and the calculation of
the Wasserstein-1 distance.
From the computational point of view, besides the totally inexpensive
projection for the gradient of basis functions $\PONE(\Trianghh)$,
finding the approximating optimal solution $(\OptTdensH,\OptPotH)$ on
surfaces requires the same effort as on the 2d Euclidean setting. A
number of improvements (multilevel approximation, hardware
acceleration, or implicit time-stepping scheme solved via
Newton-Raphson, like in~\cite{FaBe2021}) can be considered to improve the computational
efficiency of the scheme but in this paper we are more interested in
testing the accuracy of the method applied to surfaces.
In the test case where an exact solution of the \MKEQS\ is known, the
experimental convergence rate in the approximation of the solution of
the $\Lspace{1}$-\OTP\ is shown to be linear with respect to the mesh
parameter $\MeshPar$.  However, the experimental rate of convergence
in the computation of the Wasserstein-1 distance showed
superconvergence that has not yet been explained. In fact, we found an
almost cubic order of convergence and in general the discrete \DMK\
approach showed noticeably better accuracy and performance with
respect to existing algorithms.